%

\magnification=\magstep1
\input amstex
\documentstyle{amsppt}

\hsize=5.5truein
\vsize=8.5truein

\define\dt3{{\underset{\sim}\to{\delta}}^1_3}
\define\dtnew{{\underset{\sim}\to{\delta}}^1_2}
\define\dtt{{\underset{\sim}\to{\delta}}_1}

\define\om{\omega_2}
\define\ome{\omega_1}
\define\De{{\underset{\sim}\to{\Delta}}^1_2}
\define\Dee{{\underset{\sim}\to{\Delta}}_1}
\define\Dt{\Delta_1}
\define\xe{\langle x,\epsilon \rangle}
\define\aco{A\subseteq\operatorname{ORD}}
\define\Aco{\langle x,A\cap x\rangle\cong\langle\bar x, A\cap\bar x\rangle}

\define\xej{\langle\bar x_j|j<i\rangle}

\define\ilp{\{i|p(i)(\gamma )\ne(\phi ,\phi )\}}
\define\ilpp{\{i|p(i)(0)\ne(\phi ,\phi )\}}

\document

\baselineskip=.25truein

\centerline{\bf $\bold{\underset{\sim}\to{\delta}^1_2}$ Without Sharps}
\vskip10pt

{\centerline{ Sy D. Friedman,\footnote"*"{Research supported by NSF
Contract \# 9205530. Mathematical Reviews Classification Numbers: 03E15,03E35,03E55.}}}
\centerline {Department of Mathematics}
\centerline {Massachusetts Institute of Technology}
\centerline {Cambridge, MA 02139 USA}
\centerline {e-mail: sdf$\@$math.mit.edu}

\vskip2pt

\centerline{ W. Hugh Woodin, UC Berkeley}
\centerline{ Department of Mathematics}
\centerline{University of California}
\centerline{Berkeley, CA 94720 USA}
\centerline {e-mail: woodin$\@$math.berkeley.edu}

\vskip20pt

$\dtnew$ denotes the supremum of the lengths of $\De$ prewellorderings
of the reals. A result of Kunen and Martin (see Martin[77]) states that $\dtnew$ is at most $\om$ and it
is known that in the presence of sharps the assumption $\dtnew=\om$ is strong:
it implies the consistency of a strong cardinal (see Steel-Welch[?]).

In this paper we show how to obtain the consistency of $\dtnew=\om$ in the
absence of sharps, without strong assumptions.

\vskip10pt

\flushpar
{\bf Theorem.} \ {\it Assume the consistency of an inaccessible. Then it is
consistent that $\dtnew=\om$ and $\ome$ is inaccessible to reals (i\.e\.,
$\omega _1^{L[x]}$ is countable for each real $x$).}

\vskip5pt

The proof is obtained by combining the $\Dt$-coding technique of
Friedman-Velickovic [95] with the use of a product of Jensen codings of Friedman [94].

We begin with a description of the $\Dt$-coding technique.

\vskip10pt

\flushpar
{\bf Definitions.} \ Suppose $x$ is a set, $\xe$ satisfies the axiom of
extensionality and $\aco$. $x$ {\it preserves} $A$ if $\Aco$ where
$\bar x=$ transitive collapse of $x.$ For any ordinal $\delta ,$ $x[\delta
]=\{f(\gamma )|\gamma <\delta ,f\in x,f$ a function, $\gamma
\in\operatorname{Dom}(f)\}.$ $x$ {\it strongly preserves} $A$ if $x[\delta
]$ preserves $A$ for every cardinal $\delta .$ A sequence $x_0\subseteq x_1\subseteq\dots$
is {\it tight} if it is continuous and for each $i,$ $\xej$ belongs to the least $ZF^-$-model which contains  $\bar x_i$ as an element and correctly computes  card$(\bar
x_i).$

\vskip10pt

\flushpar
{\bf Condensation Condition for A.} \ Suppose $t$ is transitive, $\delta $ is
regular, $\delta \in t$ and $x\in t.$ Then:

\roster

\item"(a)" \ There exists a continuous, tight $\delta$-sequence $x_0\prec
x_1\prec\dots\prec t$ such that card$(x_i)=\delta ,$ $x \in x_0$ and $x_i$
strongly preserves $A,$ for each $i.$

\item"(b)" \ If $\delta $ is inaccessible then there exist $x_i$'s as above
but where card$(x_i)=\aleph_i.$
\endroster

The following is proved in Friedman-Velickovic [95].

\vskip5pt

\flushpar
{\bf ${\bold{\Dt}}$-Coding.} \ Suppose $V=L$ and the Condensation Condition
holds for $A.$ Then $A$ is $\Delta_1$ in a class-generic real $R,$
preserving cardinals. 

Now we are ready to begin the proof of the Theorem. Suppose $\kappa $ is
the least inaccessible and $V=L$.  Let $\langle \alpha _i|i<\kappa
^+\rangle$ be the increasing list of all $\alpha \in(\kappa ,\kappa ^+)$
such that $L_\alpha =$ Skolem hull $(\kappa )$ in $L_\alpha.$ For each
$i<\kappa ^+$ define $f_i: \kappa \longrightarrow\kappa $ by $f_i(\gamma )=$
ordertype $(ORD \cap $ Skolem hull $(\gamma )$ in $L_{\alpha _{i}}).$ By
identifying $f_i$ with its graph and using a pairing function we can think
of $f_i$ as a subset of $\kappa.$ The following is straightforward.

\vskip10pt

\flushpar
{\bf Lemma 1.} \  {\it Each $f_i$ obeys the Condensation Condition. Indeed
$\langle f_i|i<\kappa ^+\rangle$ jointly obeys the Condensation Condition
in the following sense: \ Suppose $t$ is transitive, $\delta $ is regular,
$\delta \in t,$ $x\in t.$ Then their exists a tight $\delta$-sequence
$x_0\prec x_1\prec \dots\prec t$ such that card$(x_i)=\delta ,x\in x_0$ and
each $x_i$ strongly preserves all $f_j$ for $j\in x_i$ (and if $\delta
=\kappa $ then we can alternatively require card$(x_i)=\aleph_i).$}

\vskip5pt

Now, following Friedman [94] we use a ``diagonally-supported'' product of
Jensen-style codings. For each $i<\kappa ^+$ let $\Cal{P}(i)$ be the
forcing from Friedman-Velickovic [95] to make $f_i$ $\Dt$-definable in a
class-generic real. Then $\Cal{P}$ consists of all
$p\in\prod\limits_{i<\kappa ^+}\Cal{P}(i)$ such that for infinite ordinals
$\gamma ,$ $\ilp$ has cardinality at most $\alpha $ and in addition $\ilpp$
is finite.

Now note that for successor cardinals $\gamma <\kappa $ the forcing
$\Cal{P}$ factors as $\Cal{P}_\gamma *\Cal{P}^{G_{\gamma }}$ where
$\Cal{P}_\gamma$ forces that $\Cal{P}^{G_{\gamma }}$ has the $\gamma ^+$-CC. Also
the joint Condensation Condition of Lemma 1 implies that the argument of
Theorem 3 of Friedman-Velickovic [95] can be applied here to show that
$\Cal{P}_\gamma $ is $\le\gamma$-distributive, and also that $\Cal{P}$ is
$\Delta$-distributive (if $\langle D_i|i<\kappa \rangle$ is a sequence of
predense sets then it is dense to reduce each $D_i$ below $\aleph_{i+1}$). So
$\Cal{P}$ preserves cofinalities.

Thus in a cardinal-preserving forcing extension of $L$ we have produced
$\kappa ^+$ reals $\langle R_i|i<\kappa ^+\rangle$ where $R_i$
$\Delta_1$-codes  $f_i$ and hence there are well-orderings of $\kappa $ of
any length $<\kappa ^+$ which are $\Dt$ in a real. Finally L\'evy collapse
to make $\kappa =\ome$ and we have $\dtnew=\om$, $\ome$ inaccessible to reals.
\hfill{$\dashv$}

\vskip10pt

The  above proof also shows the following, which may be of independent
interest.

\vskip10pt

\flushpar
{\bf Theorem 2.} \ Let $\dtt(\kappa )$ be the sup of the lengths of
wellorderings of $\kappa $ which are $\Dee$ over $L_\kappa [x]$ for some
$x,$ a bounded subset of $\kappa.$ Then (relative to the consistency of an
inaccessible) it is consistent that $\kappa $ be weakly inaccessible and $\dtt$
$(\kappa )=\kappa ^+.$

\vskip5pt

\flushpar
{\bf Remark.} The conclusion of Theorem 2 cannot hold in the context of
sharps: if $\kappa$ is weakly inaccessible and every bounded subset of
$\kappa$ has a sharp then $\dtt(\kappa) < \kappa^+.$ This is because
$\dtt(\kappa)$ is then the second uniform indiscernible for bounded
subsets of $\kappa$, which can be written as the direct limit of the
second uniform indiscernible for subsets of $\delta$, as $\delta$ ranges over
cardinals less than $\kappa$; so $\dtt(\kappa)$ has cardinality $\kappa$. 

\vskip10pt

Using the least inner model closed under sharp, we can also obtain the
following.

\vskip5pt

\flushpar
{\bf Theorem 3.} Assuming it is consistent for every set to have a sharp, 
then this is also consistent with $\dt3$ = $\omega_2$.

\vskip30pt

\centerline{\bf References}

\vskip10pt

\flushpar
Friedman [94] \ A Large $\Pi_{2}^1$ Set, Absolute for Set Forcings, Proceedings of the American Mathematical Society, Vol. 122, No. 1, pp. 253-256

\vskip5pt

\flushpar
Friedman-Velickovic [95] \ $\Dt$-Definability, to appear.

\vskip5pt

\flushpar
Martin [77] \ Descriptive Set Theory: Projective Sets, in {\bf Handbook of
Mathemtatical Logic}, Studies in Logic and the Foundations of Mathematics 90,
Barwise (editor), pp. 783-815. 

\vskip5pt

\flushpar
Steel-Welch [?] \ $\Sigma^1_3$ Absoluteness and the Second Uniform Indiscernible, to appear, Israel Journal of Mathematics

\vskip5pt

\enddocument